\documentstyle[amssymb,amsfonts]{amsart}

\def\Var{\hbox{\bf Var}}

\def\R{{\hbox{\bf R}}}

\def\P{{\hbox{\bf P}}}
\def\E{{\hbox{\bf E}}}

\font \roman = cmr10 at 10 true pt

\def\rank{\hbox{\rm rank}}

\def\be#1{ \begin{equation}\label{#1} }

\def\bas{\begin{align*}}
\def\eas{\end{align*}}
\def\bi{\begin{itemize}}
\def\ei{\end{itemize}}

\def\dim{{\hbox{\roman dim}}}

\def\eps{\varepsilon}
\newenvironment{proof}{\noindent {\bf Proof} }{\endprf\par}
\def \endprf{\hfill  {\vrule height6pt width6pt depth0pt}\medskip}
\def\emph#1{{\it #1}}
\def\textbf#1{{\bf #1}}

\def\BE{{\mathbf E}}

\def\BBR {{\mathbb R}}

\def\ep{{\epsilon}}
\def\hs{\hfill $\square$}

\theoremstyle{plain}
  \newtheorem{theorem}[subsection]{Theorem}
  \newtheorem{conjecture}[subsection]{Conjecture}

  \newtheorem{lemma}[subsection]{Lemma}
  \newtheorem{corollary}[subsection]{Corollary}
   
     \newtheorem{question}[subsection]{Question}

\theoremstyle{remark}
  \newtheorem{remark}[subsection]{Remark}

\theoremstyle{definition}
  \newtheorem{definition}[subsection]{Definition}

\include{psfig}

\begin{document}

\title{Random symmetric matrices are almost surely non-singular}

\author{Kevin Costello}
\address{Department of Mathematics, UCSD, La Jolla, CA 92093-0112}
\email{kcostell@@math.ucsd.edu}

\author{Terence Tao}
\address{Department of Mathematics, UCLA, Los Angeles CA 90095-1555}
\email{tao@@math.ucla.edu}
\thanks{T. Tao is supported by a grant from the Packard Foundation.}

\author{Van Vu}
\address{Department of Mathematics, UCSD, La Jolla, CA 92093-0112}
\email{vanvu@@ucsd.edu}
\thanks{V. Vu is a Sloan Fellow and is supported by an NSF Career Grant.}

\begin{abstract} Let $Q_n$ denote a random symmetric $n$ by $n$ matrix, whose upper diagonal  entries are
i.i.d. Bernoulli random variables (which take values 0 and 1 with
probability $1/2$). We prove that $Q_n$ is non-singular with probability $1-O(n^{-1/8+\delta})$ for any fixed $\delta > 0$.
The proof uses a quadratic version of
Littlewood-Offord type results concerning the concentration
functions of random variables and can be extended for more general
models of random matrices.
\end{abstract}

\maketitle

\section {Introduction}

Let $A_n$ denote a random $n$ by $n$ matrix, whose entries are
i.i.d. Bernoulli random variables, which take values 0 and $1$
with probability $1/2$. A basic  question is the following

\begin{question} \label{komlos1}  Is it true that  $A_n$ is almost surely non-singular ?
\end{question}

Here and later we say that an event holds almost surely if it
holds with probability tending to one as $n$ tends to infinity.

The above  question was answered affirmatively by Koml\'os in 1967
\cite{Kom1}. Later,  Koml\'os generalized the result  (to more
general models of random matrices)  \cite{Kom2} and also simplified
the proof \cite{Bol1}. In a recent paper \cite{TV1}, Tao and Vu
found a different proof which  leads to a sharp estimate on the
absolute value of the determinant of $A_n$.

Another popular model of random matrices is that of random symmetric
matrices; this is one of the simplest models that has non-trivial
correlations between matrix entries.  Let $Q_n$ denote a random
symmetric $n$ by $n$ matrix, whose upper diagonal entries ($q_{ij},
1\le i \le j \le n$) are i.i.d. Bernoulli random variables. It is
natural to ask

\begin{question} \label{weiss}  Is it true that  $Q_n$ is almost surely non-singular ?
\end{question}

As far as we can trace, this question was first posed by Weiss in
the early nineties. Despite its obvious similarity to Question
\ref{komlos1}, we do not know of any partial results concerning
this question,  prior to this paper.  A significant new difficulty is that the symmetry ensures that
the determinant $\det(Q_n)$ is a quadratic function of each row, as opposed to $\det(A_n)$ which is a linear function of each row.

The goal of the current  paper is to give an affirmative answer to
Question \ref{weiss}.

\begin{theorem} \label{theo:main} $Q_n$ is almost surely
non-singular. More precisely

$$ p_n := \P(Q_n \hbox{\rm \ is singular}) = O(n^{-1/8 + \delta}) $$

\noindent for any positive constant $\delta$ (the implicit constant in the $O()$ notation of course
is allowed to depend on $\delta$).

\end{theorem}

\begin{remark} \label{remark:main} The exponent $-1/8+\delta$ can be
improved somewhat by  tightening the calculation and applying  more
technical arguments. However, to improve the bound to an exponential bound (in the spirit of \cite{KKS}) seems to require new ideas; see Section \ref{section:open}.
\end{remark}

The rest of the paper is organized as follows. In the next
section, we present our approach and the key lemmas. The lemmas
will be discussed in Sections 3-5.  Section 6 is devoted to the
generalization of the result to other  models of random matrices.
We conclude by Section \ref{section:open} which contains several
open questions.

\vskip2mm

{\bf \noindent Notation.} In the whole paper, we assume that $n$
is large, whenever needed. The asymptotic notations are used under
the assumption that $n \rightarrow \infty$. $\BE$ and $\Var$
denote expectation and variance, respectively; $\log$ denotes the
logarithm with natural base.

\section{ The approach and main lemmas}

As mentioned above, there are now three different proofs of Koml\'os
1967 result on the non-singularity of $A_n$. The simpler ones are
\cite{Bol1} and \cite{TV1}. But the original (and longest) proof
from  \cite{Kom1} is what really inspires us. The key difference
between these proofs lies in the ways one generates $A_n$. In the
proofs from \cite{Bol1} and \cite{TV1} one builds up $A_n$ by
exposing the row vectors one by one and making use of the
independence of these vectors. This approach, unfortunately, is no
longer effective for  $Q_n$, as the last few rows are almost
deterministic once one has exposed all rows above them. In
\cite{Kom1}, one builds up $A_n$ by taking $A_{n-1}$ and adding a
(random) row and a (random) column. This idea turns out to be useful
for the consideration of $Q_n$. However, for $Q_n$ the additional
row and column are not independent. They are transposes of each
other and this has  become the main obstacle. We have managed to
overcome this obstacle by developing a {\it quadratic} variant of
Littlewood-Offord type results  concerning the concentration of
random variables (see Section 4).

The basic strategy is to relate the rank of $Q_n$ with the rank of
$Q_{n+1}$.  Assume that we get $Q_{n+1}$ by adding a new column
and its transpose as a new row to $Q_n$. Our starting point is the
following simple observation

\begin{equation}\label{rankn}
 \rank(Q_n) \leq \rank(Q_{n+1}) \leq \rank(Q_n) + 2.
 \end{equation}
We shall refine this by showing that if $Q_n$ is singular (so $\rank(Q_n) < n$), then $\rank(Q_{n+1})$
will equal $\rank(Q_n) + 2$ with high probability; similarly, if $Q_n$ is non-singular (so $\rank(Q_n) = n$),
then $\rank(Q_{n+1})$ will equal $\rank(Q_n)+1 = n+1$ with high probability.  These two results together will then
be easily combined with an inductive argument to show that $\rank(Q_n) = n$ with high probability.

We now turn to the details.
Let us fix a small positive constant $\ep > 0$.
We allow the implicit constants to depend on $\ep$, and we will assume that $n$ is sufficiently large
depending on $\eps$.
Set
\begin{equation}\label{N-def}
N:=n^{1-\ep}.
\end{equation}

\begin{definition}
Given $m$  vectors $\{v_1, v_2, ..., v_m\}$, a \emph{linear combination}
of the $v_i$ is a vector $v= c_1 v_1 + \dots c_m v_n$, where the
$c_i$ are real numbers. We say that a linear combination {\it
vanishes } if $v$ is the zero vector. A vanishing linear
combination has {\it degree} $k$ if exactly $k$ among the $c_i$
are non-zero.  We call a singular $n$ by $n$ matrix \textit{normal} if
its row vectors do not admit a non-trivial vanishing linear
combination with degree less than $N$. Otherwise we call the matrix
{\it abnormal}.
\end{definition}

\begin{remark} We use the terms normal and abnormal only when the
matrix in question is singular. These terms are not defined (and
we don't need them) for non-singular matrices. \end{remark}

In Section \ref{abim-sec} we shall prove that most singular matrices are normal:

\begin{lemma} \label{lemma:abnormal} The probability that $Q_n$
is singular and abnormal is $O( (2/3)^n )$.
\end{lemma}

In Section \ref{rank-sec} we shall prove

\begin{lemma} \label{lemma:rankincrease1}
Let $A$ be a (deterministic) $n$ by $n$ singular normal  matrix, and let $A'$ be
the $n+1$ by $n+1$  matrix formed by augmenting $A$ by a random
vector of length $n+1$ and its transpose. Then
$$\P (\rank (A') -\rank(A) < 2) = O(N^{-1/2}) $$
and thus
$$\P (\rank (A') -\rank(A) = 2) = 1 - O(N^{-1/2}). $$
\end{lemma}

Intuitively, these two lemmas state that in most cases, augmenting
a singular matrix by  a random vector and its transpose will
increase the rank  by exactly 2.  Note that by Bayes' identity, Lemma \ref{lemma:rankincrease1} automatically
generalizes to matrices $A$ which are random instead of deterministic, as long as the random vector which is augmenting $A$
is independent of $A$.

\vskip2mm

We now develop analogues of the above two lemmas for non-singular matrices.

\begin{definition}
A row of an $n$ by $n$ non-singular matrix is called \textit{good}
if its exclusion leads to an $(n-1) \times n$ matrix whose column
vectors admit  a nontrivial vanishing linear combination with
degree at least $N$. (In fact, there is exactly one such
combination--up to scaling--as the rank of this $(n-1) \times n$
matrix is $n-1$.) A row is {\it bad} otherwise.  We say that an $n
\times n$ non-singular matrix $A$ is \textit{perfect} if every row
in $A$ is a good row. If a non-singular matrix is not perfect, we
call it {\it imperfect}.
\end{definition}

\begin{remark} We use the terms perfect  and imperfect only when the
matrix in question is non-singular. These terms are not defined
for singular matrices. \end{remark}

In Section \ref{abim-sec} we shall prove that most non-singular matrices are perfect:

\begin{lemma} \label{lemma:imperfect} The probability that $Q_n$ is both non-singular
and imperfect is $O( (2/3)^n )$.
\end{lemma}

In Section \ref{rank-sec} we shall prove the following analogue of Lemma \ref{lemma:rankincrease1}:

\begin{lemma} \label{lemma:rankincrease2}
Let $A$ be a (deterministic) non-singular perfect symmetric matrix of size $n$,
and let $A'$ be the $(n+1) \times (n+1)$ matrix formed by augmenting
$A$ by a random (n+1)-vector of 0s and 1s, and its transpose. Then
$$\P \left(\rank (A') = \rank(A) \right) =O(N^{-1/8} ) $$
\noindent for any positive constant $\delta$, where the implicit constant can of course depend on $\delta$.  In particular, since
$$ n = \rank(A) \leq \rank(A') \leq n+1$$
we see that
$$\P \left(\rank (A') - \rank(A) = 1 \right) = 1 - O(N^{-1/8} ).$$
\end{lemma}

The last two lemmas are the non-singular counterparts of the first
two. Together, they state that if a matrix already has full rank,
augmenting it will typically produce another matrix of full rank.  Again, we can automatically generalize
Lemma \ref{lemma:rankincrease2} to the case when $A$ is random and independent of the augmenting row.

Let us assume these lemmas for the moment and conclude the proof of Theorem \ref{theo:main}.

Consider a random matrix $Q_n$. We embed it into a sequence $\{Q_1, Q_ 2, ...\}$ of random
matrices, where $Q_{n+1}$ is formed from $Q_n$ by adding a random
vector of 0s and 1s (independent of $Q_n$) of length $n+1$  as the last column, and its transpose as
the last row.

Define the (somewhat artificial) random variable $X_n$ by setting $X_n = 0$ if $Q_n$ is
non-singular (thus $\rank(Q_n) = n$), and $X_n = (1.1)^{n -
\rank(Q_n)}$ otherwise.  Thus $X_n$ ranges between 0 and
$(1.1)^n$.  We have the following decay estimate for the
expectation $\E(X_n)$ of $X_n$.

\begin{lemma}   $\E(X_{n+1}) \leq 0.99 \E(X_n) + O( N^{-1/8} )$.
\end{lemma}

\begin{proof}  For any $0 \leq j \leq n$, let $A_j$ be the event that $Q_n$ has rank $n-j$, and that $A_n$ is neither abnormal (if $j>0$) nor
imperfect (if $j=0$).  By Bayes' identity and Lemmas \ref{lemma:abnormal}, \ref{lemma:imperfect}, we have
$$ \E(X_n) = \sum_{j=0}^n (1.1)^j \P( A_j) + O( (1.1)^n (2/3)^n )$$
and
$$ \E( X_{n+1} ) = \sum_{j=0}^n \E( X_{n+1} | A_j ) \P( A_j ) + O( (1.1)^n (2/3)^n ).$$
Now let us condition on the event $A_0$, thus $Q_n$ is
non-singular.  From Lemma \ref{lemma:imperfect} we see that
$Q_{n+1}$ has rank $n$ with probability $O(N^{-1/8})$, and rank
$n+1$ otherwise.  Thus
$$ \E( X_{n+1} | A_0 ) = O( N^{-1/8} ).$$
Now let $1 \leq j \leq n$ and condition on the event $A_j$, thus
$Q_n$ is singular with rank $n-j$.  From Lemma
\ref{lemma:abnormal} and Lemma \ref{lemma:rankincrease1} we see
that $Q_{n+1}$ has rank $n-j+2$ with probability $1 -
O(N^{-1/2})$, and has rank $n-j$ or $n-j+1$ otherwise.  Thus
\begin{align*}
\E( X_{n+1} | A_j ) &\leq \E( 2^{n+1-\rank(Q_{n+1})} | A_j ) \\
&\leq (1.1)^{j-1} + O( N^{-1/2} ) (1.1)^{j+1}\\
&\leq 0.99 (1.1)^j
\end{align*}
if $N =n^{1-\ep}$ is large enough.  Putting all these estimates
together, and noting that $$O( (1.1)^n (2/3)^n ) = O( N^{-1/8} ),
$$ we obtain the claim.
\end{proof}

From the above lemma and an easy induction, we see that
$$ \E(X_n) = O( N^{-1/8} )$$
for all large $n$.  From Markov's inequality we then see that
$$ \P( Q_n \hbox{ singular} ) = \P(X_n \geq 1) =  O( N^{-1/8} ).$$
Theorem \ref{theo:main} then follows from the definition of $N$.

It remains to prove Lemmas
\ref{lemma:abnormal}--\ref{lemma:rankincrease2}.  This will be done in the next few sections.
Among these lemmas, the first three are variants of
lemmas from \cite{Kom1} and are relatively simple. The proof of
Lemma \ref{lemma:rankincrease2} is a somewhat more complicated and
requires a new machinery, discussed in Section \ref{section:LO}.

\section{Proof of Lemmas \ref{lemma:abnormal} and \ref{lemma:imperfect}  }\label{abim-sec}

The two proofs are similar and rely on the following simple
observation from \cite{Kom1} (which has also been used in
\cite{KKS}, \cite{TV1}, \cite{TV2}):

\begin{lemma} \label{lemma:subspace}  Let  $H$ be a linear subspace in $\BBR^n$ of
dimension at most $d \le n$. Then it contains at most $2^d$
vectors from $\{0,1\}^n$.
\end{lemma}

\begin{proof} The space $H$ is spanned by the row vectors of a  $d' \times n$ full-ranked matrix,
where $d' := \dim(H) \leq d$.  This matrix has at least one
non-singular $d' \times d'$ minor, thus there exists a set of $d'$
co-ordinates of $\BBR^n$ which can be used to parameterize $H$.
But in $\{0,1\}^n$, these co-ordinates take only $2^{d'} \leq 2^d$
values, and the claim follows.
\end{proof}

\vskip2mm

\begin{proof}[of Lemma \ref{lemma:abnormal}] For any $1 \leq d \leq n$. Let $g(n,d)$ be the probability that the row vectors of
$Q_n$ admit a nontrivial vanishing linear combination of degree
$d$.  For $d=1$ we have the easy bound $g(n,1) \leq n 2^{-n}$, since $g(n,d)$ is simply the probability that one of the rows of $Q_n$ is
entirely zero.  Now take $d \geq 2$.  To bound $g(n,d)$ from above, notice that by symmetry and the union bound
we have the crude estimate
$$g(n,d) \le {n \choose d} h(n,d) \leq n^{d} h(n,d) $$

\noindent where $h(n,d)$ is the probability that the first $d$
rows admit such combination. This means that if we fix the first
$d-1$ row vectors, then the $d^{th}$ row vector lies in the
subspace spanned by these vectors.  The same claim is true if we
delete the first $d-1$ columns from $Q_n$ (we need to do this as
$Q_n$ is symmetric). The remaining entries in the $d^{th}$ row
vector are now distributed independently in $\{0,1\}^{n-d-1}$, and
so by Lemma \ref{lemma:subspace} the probability of lying in the
span of the first $d-1$ row vectors is at most
 is at most $2^{d-1} / 2^{n-d-1}$. Thus the probability that $Q_n$ is singular and abnormal is at most
$$\sum_{d=1}^N g(n,d) \le n 2^{-n} + \sum_{d=2}^N n^{d+1} 2^{d-1} / 2^{n-d-1} = O( (2/3)^n )$$
by the definition \eqref{N-def} of $N$.
\end{proof}

\vskip2mm

\begin{proof}[of Lemma \ref{lemma:imperfect}] Let $b(n)$
be the probability that the last row of $Q_n $ is bad. By symmetry
and the union bound, the probability that $Q_n$ is non-singular and
imperfect is at most $nb(n)$. We can bound $b(n)$ using the same
argument as in the previous proof, with a slight modification; the
column vectors have length $n-1$ so we need to replace $n$ by $n-1$,
but this does not affect the bound. We omit the details.
\end{proof}

\section{A quadratic Littlewood-Offord inequality}\label{section:LO}

Let us start by the  following classical result, proved by
Erd\"os, which strengthens an earlier result of Littlewood and
Offord.

\begin{theorem}[Linear Littlewood-Offord inequality]\label{theo:LO} \cite{erdos}
Let $z_1,\ldots,z_n$ be i.i.d. random variables which take values
$0$ and $1$ with probability $1/2$. Let $a_1,\ldots,a_n$ be real
deterministic coefficients, with $|a_i| \geq 1$ for at least $k
\ge 1$ values of $i$.  Then for any interval $I \subset \R$ of
length 1, we have
$$\P\left(\sum_{i=1}^n a_i z_i \in I\right)  = O\left( k^{-1/2}\right) $$
where the implied constant is absolute.
\end{theorem}

Roughly speaking, the theorem  says that linear random sums cannot
concentrate on small intervals if the coefficients of the
underlying linear form are large.

\begin{remark} \label{remark:LO}
There are a number of far reaching generalizations and interesting
refinements of Theorem \ref{theo:LO} (see e.g. \cite{Hal} and the
references therein).  We mention some rather trivial ones here (which we will need later).  Firstly we can replace
the unit interval $I$ by any other interval of length $O(1)$ (at the cost of changing the implied constant
in $O((1+k)^{-1/2})$, of course), by covering such an interval by unit intervals.  Similarly, we may scale
the constraint $|a_i| \geq 1$ and replace it by $|a_i| \geq c$ for some other $c > 0$, again at the cost of letting
the implied constant depend on $c$.  Finally, one can replace the distribution of the $z_i$ with
the distribution $\P(z_i=0) = \alpha, \P(z_i=1) =\beta, \P(z_i=-1) =
\gamma$, where $\alpha, \beta,\gamma$ are  non negative  constants
summing up to one and $\alpha <1$. The implied constant will then of course
 depend on $\alpha, \beta, \gamma$.
\end{remark}

To conclude the proof of Theorem \ref{theo:main}, we  need to
generalize Theorem \ref{theo:LO} in a direction  different from
what has been done before. Instead of considering a linear form,
we are going to consider a \emph{quadratic} form of the $z_i$. (In fact,
our method works for polynomials of any fixed degree, by iterating the argument below.) Consider
random variables $z_i$ as in Theorem \ref{theo:LO} and define
the quadratic random variable
\begin{equation}\label{QDEF}
Q := \sum_{1\le i, j \le n} c_{ij} z_i z_j.
\end{equation}

The main result of this section is the following quadratic generalization of Theorem \ref{theo:LO}.

\begin{theorem}[Quadratic Littlewood-Offord inequality]\label{theo:quadLO}
Let the quadratic random variable $Q$ be as in \eqref{QDEF}, let
$\{1, \dots, n\} = U_1 \cup U_2$ be any non-trivial partition, and
let $S$ be any non-empty subset of $U_1$. For each $i \in S$, let
$d_i$ be the number of indices $j \in U_2$ such that $|c_{ij}|
\geq 1$. Suppose that $d_i \ge 1$ for each $i \in S$. Then for any
interval $I$ of length 1, we have
$$\P(Q \in I) = O\left( |S|^{-1/2} + |S|^{-1} \sum_{i \in S} d_i^{-1/2} \right)^{1/4},$$
The implied constant is absolute.
\end{theorem}

It is unlikely that the bound on the right-hand side is best possible, but for us, any bound which decays to zero when the number of
large coefficients $c_{ij}$ goes to infinity will suffice.

The proof of Theorem \ref{theo:quadLO} is lengthy and will be given later.  Assuming it for the moment, we have the following
corollary:

\begin{corollary} \label{cor:quadLO1} Let $Q$ be as in \eqref{QDEF}, and suppose that there is a set $U \subseteq \{1,\ldots,n\}$
of cardinality $|U| \geq m \geq 2$ such that for each
$i \in U$, there are $m$ indices $j \in \{1,\ldots,n\}$ where
$|c_{ij}| \ge 1$. Then for any interval $I$ of length 1
$$\P(Q \in I)  \le O(m^{-1/8}).$$
The implied constant is absolute.
\end{corollary}

\begin{proof} Without loss of generality we may take $m$ to be even.  Let $U_1$ be an arbitrary
subset of $U$ of cardinality $m/2$ and write $U_2 := \{1,\ldots,n\} \backslash U_1$,
then for any $i \in U_1$ there exists at least $m/2$
indices $j \in U_2$ for which $|c_{ij}| \geq 1$.  Applying Theorem \ref{theo:quadLO} with $S := U_1$, we conclude
\begin{align*}
\P(Q \in I) &= O( (m/2)^{-1/2} + (m/2)^{-1} \sum_{i \in U_1} (m/2)^{-1/2} )^{1/4} \\
 &=O(m^{-1/8})
\end{align*}
as desired.
\end{proof}

\noindent By rescaling the above corollary, we obtain the following discrete version.

\begin{corollary} \label{cor:quadLO2}
Let $Q$ be as in \eqref{QDEF}, and suppose that there are at least $m$ indices  $i$ such that for each
$i$ there are $m$ indices $j$ where $|c_{ij}| \neq 0$.
Then
$$\P(Q = 0) = O(m^{-1/8})$$
where the implied constant is absolute.
\end{corollary}

This Corollary will be the one we use to establish Lemma
\ref{lemma:rankincrease2}.

\subsection{Proof of Theorem \ref{theo:quadLO}}

We now prove Theorem \ref{theo:quadLO}.  As a first attempt to prove this theorem, one might
try to view the quadratic form $Q$ as a linear form
\begin{equation}\label{Q-factor}
Q = \sum_{i=1}^n Q_i z_i,
\end{equation}
where the coefficients $Q_i$ are themselves linear form random
variables $Q_i := \sum_{j=1}^n c_{ij} z_j$.  Thus one might hope
to obtain Theorem \ref{theo:quadLO} from two applications of
Theorem \ref{theo:LO}.  Unfortunately, there is a serious
obstruction to this strategy, because the coefficients
$Q_1,\ldots,Q_n$ are not independent of the variables
$z_1,\ldots,z_n$.  However, we can get around this obstacle by the
following decoupling lemma, which relies on the Cauchy-Schwarz
inequality.

\begin{lemma}[Decoupling lemma]\label{lemma:decoupling} Let $X$ and $Y$ be
 random variables and $E=E(X,Y)$ be an event depending on $X$ and
$Y$. Then

$$\P( E(X,Y)) \le\P( E(X,Y) \wedge E(X',Y) \wedge E(X,Y')\wedge
E(X',Y'))^{1/4} $$

\noindent where $X'$ and $Y'$ are independent copies of $X$ and $Y$,
respectively.  Here we use $A \wedge B$ to denote the event that $A$ and $B$ both hold.
\end{lemma}

\begin{remark}  This lemma is a probabilistic analogue of the well-known result in extremal graph theory,
that if a bipartite graph connecting $n$ and $m$ vertices contains at least $cnm$ edges for some $0 \leq c \leq 1$,
then it also contains at least $c^4 n^2 m^2$ copies of the four-cycle $C_4$, where we include degenerate four-cycles.
Indeed, the two results are easily shown to be equivalent.  This decoupling lemma also plays
the role of the van der Corput lemma used
in Weyl's estimation of exponential sums with quadratic (or more generally
polynomial) phases; indeed it is quite likely that one could obtain
an estimate very similar to Theorem \ref{theo:quadLO} by means of these
techniques (combined with Ess\'een's concentration inequality), however we have
chosen a more elementary combinatorial approach here.
\end{remark}

\begin{proof} Let us first consider the case when $X$ takes a finite
number of values $x_1, \dots, x_n$ and $Y$ takes a finite number of
values $y_1, \dots, y_m$.  From Bayes' identity we have
$$ \P(E(X,Y)) = \sum_{i=1}^n \P(E(x_i,Y)) \P(X=x_i)$$
and
$$ \P(E(X,Y) \wedge E(X,Y')) = \sum_{i=1}^n \P(E(x_i,Y))^2 \P(X=x_i)$$
and hence by the Cauchy-Schwarz inequality
$$ \P(E(X,Y)) \leq \P(E(X,Y) \wedge E(X,Y'))^{1/2}.$$
Similarly, we have
$$ \P(E(X,Y) \wedge E(X,Y')) = \sum_{j=1}^m \sum_{j'=1}^m \P( E(X,y_j) \wedge E(X,y_{j'}) ) \P( Y = y_j ) \P( Y = y_{j'} )$$
and

\begin{align*}  & \,\, \P(E(X,Y) \wedge E(X,Y') \wedge E(X',Y) \wedge E(X',Y')) \\ &=
\sum_{j=1}^m \sum_{j'=1}^m \P( E(X,y_j) \wedge E(X,y_{j'}) )^2 \P(
Y = y_j ) \P( Y = y_{j'} )\end{align*}  so by Cauchy-Schwarz again
$$ \P(E(X,Y) \wedge E(X,Y')) \leq  \P(E(X,Y) \wedge E(X,Y') \wedge E(X',Y) \wedge E(X',Y'))^{1/2}.$$
Combining these two applications of Cauchy-Schwarz, we obtain the claim.  The general case when $X$ and $Y$ could be take a countable or
uncountable number of values then follows, either by a discretization argument, or by replacing the sums with integrals and using Fubini's theorem;
we omit the details, since for our application we only need the case when $X,Y$ take finitely many values.
\end{proof}

We return to the task of proving Theorem \ref{theo:quadLO}. Let $Z
\in \{0,1\}^n$ be the random variable $(z_1, \dots, z_n)$.
Consider the quadratic form $Q(Z)= Q(z_1, \dots, z_n)$ defined by
\eqref{QDEF}, and fix a non-trivial partition $\{1, \dots, n\} =
U_1 \cup U_2$ and a non-empty subset $S$ of $ U_1$. Let $I$ be an
interval of length 1. We need to prove that
$$\P( Q(Z) \in I)^4  = O( |S|^{-1/2} + |S|^{-1} \sum_{i \in S} d_i^{-1/2} ).$$
 Define $X :=(z_i)_{i \in U_1}$ and $Y
:=(z_i)_{i \in U_2}$. We can write $Q(Z)=Q(X,Y)$. Let
$z_i'$ be an independent copy of $z_i$ and set $X' :=(z'_i)_{i \in U_1}$
 and $Y' := (z'_i)_{i \in U_2})$.  Applying Lemma \ref{lemma:decoupling}, we see that it suffices to show that
$$ \P( Q(X,Y), Q(X,Y'), Q(X',Y), Q(X',Y') \in I ) =
O( |S|^{-1/2} + |S|^{-1} \sum_{i \in S} d_i^{-1/2} ).$$
A simple   calculation shows that the random variable
$$ R := Q(X,Y) - Q(X',Y)- Q(X,Y') + Q(X',Y') $$
can be written as
\begin{align*}
R &= \sum_{i \in U_1} \sum_{j \in U_2}
 c_{ij} (z_i-z_i') (z_j-z_j') \\
 &= \sum_{i \in U_1} R_i w_i
\end{align*}
where for $i \in U_1$, $w_i$ is the random variable $w_i := z_i - z'_i$, and $R_i$ is the random variable
$$ R_i := \sum_{j \in U_2} c_{ij} w_j.$$
We have eliminated the coupling problem in the factorization \eqref{Q-factor}, because the
random variables $(R_i)_{i \in U_1}$ are independent of the random variables $(w_i)_{i \in U_1}$.

Consider the four events $Q(X,Y) \in I, Q(X',Y) \in I, Q(X,Y') \in
I$ and $Q(X',Y') \in I$. If all of these hold, then $R$ lies in
the interval $J := 2I -2I$ of length 4. Thus, it suffices to show
that
$$ \P( R \in J ) = O( |S|^{-1/2} + |S|^{-1} \sum_{i \in S}  d_i^{-1/2} ).$$

Recall that for each $i \in U_1$,
$d_i$ be the number of coefficients $j \in U_1$ for which
$|c_{ij}| \geq 1$.  For each $i\in S \subseteq U_1$, we may apply Theorem \ref{theo:LO} (and Remark \ref{remark:LO}) to the
random variable $R_i$ to obtain
\begin{equation}\label{rgi}
\P( |R_i| < 1 ) = O( d_i^{-1/2} ).
\end{equation}
By the union bound we thus have the crude estimate
$$ \P( |R_i| \geq 1 \hbox{ for all } i \in S ) = 1 -  O( \sum_{i \in S} d_i^{-1/2} ).$$
This use of the union is somewhat wasteful and we can do better by
invoking the second moment method.   For each $i \in S$, let $I_i$ be the indicator
variable of the event $|R_i| \geq 1$, thus $I_i=1$ when $|R_i| \geq 1$ and $I_i = 0$ otherwise.
Thus \eqref{rgi} can be rewritten as
$$ \E(I_i) = 1 - O( d_i^{-1/2} )$$
and hence by linearity of expectation
$$ \E( \sum_{i \in S} I_i ) = |S| - O( \sum_{i \in S} d_i^{-1/2} ).$$
Also, since $d_i \geq 1$, we have at least one $j \in U_2$ for which $|c_{ij}| \geq 1$, which
easily implies that $\E(I_i) \geq 1/2$.  Thus we also have
$$ \E( \sum_{i \in S} I_i ) \geq |S|/2.$$
Next we compute the variance of $\sum_{i \in S} I_i$:
\begin{align*}
\Var (\sum_{i \in S} I_i) &=
\E((\sum_{i \in S} I_i)^2) - \E(\sum_{i \in S} I_i)^2 \\
&\leq |S|^2 - (|S|- O\left(\sum_{i\in S} d_i^{-1/2})\right)^2 \\
&= O(  |S| \sum_{i \in S} d_i^{-1/2}).
\end{align*}
\noindent By Chebyshev's inequality, we conclude
$$\P( \sum_{i \in S} I_i \le \frac{1}{4}|S| ) \le \frac{ 16 \Var (\sum_{i \in S} I_i)}{ |S|^2}
=O( \frac{1}{|S|} \sum_{i \in S} d_i^{-1/2}). $$ Thus with
probability $1 - O( \frac{1}{|S|} \sum_{i \in S} d_i^{-1/2})$, we
have $|R_i| \geq 1$ for at least $|S|/4$ values of $i \in S$.

Let us now temporarily condition the $R_i$ to be fixed for all $i
\in U_1$, and assume that $|R_i| \geq 1$ for at least $|S|/4$
values of $i \in S$.  Applying Theorem \ref{theo:LO} (and Remark
\ref{remark:LO}) to $R = \sum_{i \in U_1} w_i R_i$ (treating the
$R_i$ as fixed coefficients), we have the conditional probability
estimate
$$ \P( R \in J | R_i \hbox{ fixed } ) = O( |S|^{-1/2} ).$$
By the preceding discussion and Bayes identity, we thus have
$$ \P( R \in J ) = O( |S|^{-1/2} ) + O( \frac{1}{|S|} \sum_{i \in S} d_i^{-1/2})$$
as desired.
\endprf

\section{Proof of Lemmas  \ref{lemma:rankincrease1}
and \ref{lemma:rankincrease2} }\label{rank-sec}

{\bf \noindent Proof of Lemma \ref{lemma:rankincrease1}.} Let $A$
be a normal symmetric singular $n \times n$ matrix of rank $d$.
Let $v_i$ be the $i$th row vector of $A$. Without loss of
generality, we can assume that $v_1, \dots, v_d$ are linearly
independent. Thus, the last row vector $v_n$ can be written as a
linear combination of these vectors in a unique way

$$v_n = \sum_{i=1}^d c_i v_i. $$

\noindent As $A$ is normal,  by definition at least $N$ among the
coefficients $c_i$ are non-zero.

Consider the addition of a random $(0,1)$ column  of length $n$ to
$A$. Each of the row vectors $v_i$ receives a new (random)
coordinate and becomes a new vector $v_i'$. Clearly, $v_1', \dots,
v_d'$ are still independent. If the new matrix $A'$ fails to have
a larger rank than $A$, then the last row $v_n'$ must remain
within the span of $v_1', \dots, v_d'$. By considering the first
$n$ coordinates,  the only way this can happen is if

$$v_n'=\sum_{i=1}^d c_j v_j'. $$

\noindent This implies that the last coordinate  $x_{n+1}$ of
$v'_n$ satisfies

\begin{equation} \label{equa:linearbound} x_{n+1} =\sum_{i=1}^d c_i y^i_{n+1} ,
\end{equation}

\noindent where $y^i_{n+1}$ is the last coordinate of $v_i'$.
Since $x_{n+1}$ and  $y^i_{n+1}$ are i.i.d $(0,1)$ random
variables and at least $N$ of the $c_i$ are non-zero, Theorem
\ref{theo:LO} (see also Remark \ref{remark:LO}) implies that the
probability that \eqref{equa:linearbound} holds is $O(N^{-1/2})$.
Thus, we can conclude that with probability $1-O(N^{-1/2})$, the
new column increases the rank by one. If adding the new column
increases the rank by one, then by the fact that  $A$ is symmetric,
adding the column and its transpose as a new row increases the
rank of $A$ by 2 (regardless the value of the last diagonal
entry), concluding the proof. \hs

\vskip2mm

{\bf \noindent Proof of Lemma \ref{lemma:rankincrease2}.} Let $A$
be a perfect non-singular symmetric  matrix of order $n$. Let
$A'$ be the $n+1$ be $n+1$ symmetric  matrix obtained from $A$ by
adding a new random $(0,1)$ column $u$ of length $n+1$ as the
$n+1$st column and its transpose as the $n+1$st row.

Let $x_1, \dots, x_{n+1}$ be the coordinates of $u$; $x_{n+1}$ is
the low-right diagonal entry of $A'$. The determinant for $A'$ can
be expressed as

$$(\det A) x_{n+1} +\sum_{i=1}^n c_{ij} x_i x_j $$

\noindent where $c_{ij}$ is the $ij$ cofactor of $A$. We can
rewrite $\det A'$ as

$$Q(x_1, \dots, x_{n+1}) = (\det A) x^2_{n+1} +\sum_{i=1 ^n } c_{ij} x_i x_j $$

\noindent thanks to the fact that $x_{n+1} ^2= x_{n+1} $. We are
going to bound the probability that $Q=0$.

In order to apply Corollary \ref{cor:quadLO2}, we next show that
for each $1 \le i \le n$, many among the $c_{ij}$ are not zero.

 Since $A$ is non-singular, dropping  the $i$th row (for any $1\le i \le n$) results in  an  $n-1 \times n$
 matrix whose columns
 admit a unique (up to scaling) vanishing linear combination $\sum_{j=1}^n a_j
 u_j$.
 As $A$ is perfect,  at least $N$  among the coefficients $a_j$ are non-zero. For each
$j$ where $a_j \neq 0$, dropping both the $i$th row and the $j$th
column must result in a full rank matrix of order $n-1$. Thus
$c_{ij} \neq 0$. Thus, we can conclude  that for each $1\le i \le
n$, there are at least $N$ indices $j$ where $c_{ij} \neq 0$. The
claim of the lemma follows by applying Corollary \ref{cor:quadLO2}
with $m=N$. \hs

\section{More general results}\label{gen-sec}

In this section we briefly discuss (without detailed proofs)
several easy extensions of the method to yield some variants and generalizations of
our results.

\subsection{Generalizations of Theorem \ref{theo:quadLO}}

Theorem \ref{theo:quadLO} and Corollary \ref{cor:quadLO1} can be
extended to polynomials with arbitrary degree. One such extension
reads as follows:

\begin{theorem}\label{theo:polyLO}
Let $z_1, z_2, \dots z_n $be i.i.d. random variables which take
values 0 and 1 with probability $1/2$. Let $k$ be a fixed positive
integer. Let

\begin{equation*}
f:=\sum_{1 \le i_1 , i_2 , \dots i_k \le n} c_{i_1, i_2, \dots
i_k} z_{i_1}z_{i_2} \dots z_{i_k}.
\end{equation*}

where  at least $ n^{k-1}m$ of the coefficients $c_{i_1, i_2,
\dots i_k}$ are at least 1 in absolute value. Then for any
interval $I$ of length 1

\begin{equation*}
{\hbox{\bf P}}(f \in I) = O(\frac{1}{m^{a_k}} ).
\end{equation*}

where $a_k= 2^{-(k^2+k)/2}$ and the implicit constant in $O$
depends on $k$.

\end{theorem}

The proof proceeds via induction on $k$, with the base case being
the classical Littlewood-Offord lemma and the inductive step closely
following that of Theorem \ref{theo:quadLO}, including the use of the
following generalization of the decoupling lemma (also proven by
induction on $k$):

\begin{lemma}[Decoupling lemma]\label{lemma:decouplinggeneral} Let $X_1, \dots, X_k$  be
 random variables and $E=E(X_1, \dots, X_k)$ be an event depending on the $X_i$. Then

$$\P( E(X_1, \dots, X_k) ) \le\P( \bigwedge_{S \subset \{1, \dots, k\}  } E(X_1^S, \dots, X_k^S))^{1/2^k} $$

\noindent where $X_i^S := X_i $ if $i \in S$ and $X_i^S := X_i'$, an
independent copy of $X_i$ if $i \notin S$.
\end{lemma}

Theorem \ref{theo:quadLO} can also be extended to more general
classes of variables than the Bernoulli random variable (taking values $0$ and $1$
with equal probability) by a nearly identical proof, with the main
difference being that the base case, Theorem \ref{theo:LO}, must be replaced by
\cite[Theorem 4]{Hal}.

\subsection{Generalizations of Theorem \ref{theo:main}}

We say that a random variable $\xi$ has the $\rho$-property if

$$\max_{c \in \BBR} \P (\xi=c) \le \rho. $$

Let $\xi_{ij}$, $1 \le i, j \le n$ be independent random variables.
Assume that there is a constant $\rho<1$ (not depending on $n$) such
that for all $1 \le i < j \le n$, $\xi_{ij}$ has the
$\rho$-property. Observe  that we do not require $\xi_{ij}$ be
identical, and that furthermore we do not place any requirements on
the diagonal elements of the matrix.

\begin{theorem} \label{theo:gen} Let $\xi_{ij}$, $1 \le i \le j \le n$
be as above. Let $Q_n$ be the random symmetric matrix with upper
diagonal entries $\xi_{ij}$. Then $Q_n$ is non-singular with
probability $1-O(n^{-1/8 + \delta})$, where the implicit constant
depends only on $\rho$ and $\delta$.
\end{theorem}

To prove this result, it suffices to show that analogues of Lemmas
\ref{lemma:abnormal}--\ref{lemma:rankincrease2} still hold for this
more generalized model.  Lemmas \ref{lemma:abnormal} and
\ref{lemma:imperfect} (with 2/3 replaced by any $\delta$ with
$\rho<\delta<1$) follow from the same argument as in the original
theorem, except that Lemma \ref{lemma:subspace} must be replaced by

\begin{lemma}
Let $H$ be a linear subspace in $\BBR^n$ of dimension at most $d \le
n$.  Let $v$ be a vector whose entries are independent random
variables all but one of which have the $\rho$ property.  Then

$$\P (v \in H)< \rho^{n-d-1}$$
\end{lemma}
\begin{proof}
As before, $H$ can be parameterized by some set of $d' \le d$
coordinates.  Once those coordinates of $v$ are known, the remaining
coordinates can each take on at most one value for all $(v \in H)$,
giving a necessary set of $(n-d')$ independent events, $(n-d'-1)$ of
which have probability at most $\rho$.
\end{proof}

The proof of Lemma \ref{lemma:rankincrease1} also goes through,
except that Theorem \ref{theo:LO} must be replaced by the following
rescaled version of the $d=1$ case of \cite[Theorem 4]{Hal}:

\begin{lemma}
Let $z_1,\ldots,z_n$ be independent random variables with the
$\rho$ property. Let $a_1,\ldots,a_n$ be real deterministic
coefficients, with $a_i \neq 0$ for at least $k$ values of $i$.
Then for any interval $c \in \R$, we have
$$\P(\sum_{i=1}^n a_i z_i =c)  = O\left((1+k)^{-1/2}\right) $$
where the implied constant depends only on $\rho$.
\end{lemma}

A nearly identical decoupling argument now proves an analogue of
Theorem \ref{theo:quadLO}, with $d_i$ now taken for each $i$ to be
the number of $j$ for which $c_{ij}$ is nonzero.  Corollary
\ref{cor:quadLO2} (with the implied constant now depending only on
$\rho$) and Lemma \ref{lemma:rankincrease2} now follow as before.

\section{Open questions} \label{section:open}

Let us conclude this section with a few open questions. From a
quantitative point of view, there are two natural ways to
strengthen both Questions \ref{komlos1} and \ref{weiss}.

\begin{question} \label{tv1}  Give an estimate for the
determinant.
\end{question}

\begin{question} \label{tv2}  Give an estimate for the
probability that the matrix is singular.
\end{question}

In fact, Question \ref{tv1} seems to be  the motivation of
Koml\'os for his original paper \cite{Kom1} (see the title of
that paper)  which started this line of research. There are
several partial results concerning the model $A_n$. In the rest of
this section, it is more convenient to assume that the entries of
$A_n$ (and $Q_n$) take value $1$ and $-1$ (rather than $1$ and
$0$). Under this condition,  Tao and Vu  \cite{TV1} showed that
almost surely $\det A_n$ has absolute value $n^{(1/2-o(1))n}$. We
conjecture that a similar bound holds for $Q_n$.

\begin{conjecture} \label{determinant} Almost surely, $|\det Q_n| = n^{(1/2-o(1))n}$.
\end{conjecture}

Regarding Question \ref{tv2}, Kahn, Koml\'os and Szemer\'edi
\cite{KKS} proved that the singular probability of $A_n$ is
$O(.999^n)$. This bound has recently been improved \cite{TV2} to
$(3/4+o(1))^n$. The conjectured bound is $(1/2+o(1))^n$. We
conjecture that the same bound holds for $Q_n$.

\begin{conjecture} \label{determinant-2} The probability that  $Q_n$
is singular is $(1/2+o(1))^n$.
\end{conjecture}

By considering the probability that the first two rows are equal,
it is easy to see that $(1/2+o(1))^n$ is a lower bound (one can
actually makes a more precise conjecture similar to the case with
$A_n$). The proof in this paper showed a upper bound $O(n^{-c})$
for some positive constant $c$.

The main obstacle in these questions is the fact that the row
vectors of $Q_n$, unlike those of $A_n$, are not independent. In
fact, if one exposes these vectors one by one, then the last few
vectors are almost deterministic. The independence among the row
vectors are critical in all previous papers \cite{KKS, TV1, TV2}.
It so  seems to require a new idea to attack these conjectures.

\vskip2mm

{\it Acknowledgement.} We would like to thank G. Kalai for
communicating the problem.

\end{document}